\documentclass{amsart}
\usepackage{latexsym,amsmath,amssymb}

\theoremstyle{definition}
\theoremstyle{remark}

\numberwithin{equation}{section}

\begin{document}

\title[SUBSTITUTION OPERATORS ON THE SPACE OF FORMAL POWER SERIES]
{SUBSTITUTION OPERATORS ON THE SPACE OF FORMAL POWER SERIES WITH A
WEIGHTED CAUCHY PRODUCT }
\author{\bf Y. Estaremi  and M. R. Jabbarzadeh }
\address{ Y. Estaremi and M. R. Jabbarzadeh }
\address{Department of Mathematics, Payame Noor University, P. O. Box 19395-3697, Tehran, Iran.}

\address{ Faculty of Mathematical Sciences, University of Tabriz, P. O. Box: 5166615648,
Tabriz, Iran}
 \email{yestaremi@pnu.ac.ir} \email{mjabbar@tabrizu.ac.ir }

\thanks{}



\date{}

\dedicatory{}

\commby{}

\begin{abstract}

Boundedness of $\diamond$-substitution operator $u\diamond
C_\varphi$ on $\ell^{p}(\beta)$, with a weighted Cauchy product
$\diamond$, is investigated by an inductive argument on the pair
$(u, \varphi)$.
\end{abstract}

\maketitle

\section {\sc\bf Introduction}

 \vspace*{0.3cm}Let
$\{\beta_{n}\}^{\infty}_{n=0}$ be a sequence of positive numbers
with $\beta(0)=1$ and $1\leq p <\infty$. We consider the space of
sequences $f=\{\hat{f}(n)\}$ such that
$$\|f\|_{\beta}^p=\sum^{\infty}_{n=0}|
\hat{f}(n)|^{p}\beta(n)^{p}<\infty.$$ We shall use the formal
notation $f(z)= \sum^{\infty}_{n=0}\hat{f}(n)z^{n}$ whether or not
the series converges for any complex values of $z$. These are called
formal power series.
 Throughout this paper, we consider the space
$\ell^{p}(\beta)$ to be defined by
$$\ell^{p}(\beta)=\{f: f(z)=\sum^{\infty}_{n=0}\hat{f}(n)z^{n}, \
\|f\|_{\beta}<\infty\} .$$ Note that if
$\lim_{n}{\beta(n+1)}/{\beta(n)}=1$ or
$\liminf_{n}\beta(n)^{{1}/{n}}=1$, then $\ell^{p}(\beta)$
consists of functions analytic on the open unit disc
$\mathbb{D}=\{z\in\mathbb{C}: |z|<1\}$. The Hardy, Bergman and
Dirichlet spaces can be viewed in this way when $p = 2$ and
respectively $\beta(n) = 1$, $\beta(n)= (n + 1)^{-1/2}$ and
$\beta(n) = (n + 1)^{1/2}$.  For a beautiful exposition of the
study of classic $\ell^{p}(\beta)$ spaces, see \cite{c, sh, z} and
the references therein.

\vspace*{0.3cm} Let $\{\delta_{n}\}^{\infty}_{n=0}$ be
  a sequence of positive numbers with $\delta_{0}=1.$
Given arbitrary two functions
$f(z)=\sum^{\infty}_{n=0}\hat{f}(n)z^{n}$
  and $g(z)=\sum^{\infty}_{n=0}\hat{g}(n)z^{n}$ of the space $\ell^{p}(\beta)$,
  define the following generalized  Cauchy product series

 \begin{eqnarray}
  \ (f\diamond \
  g)(z)=\sum^{\infty}_{n=0}\sum^{\infty}_{m=0}\frac{\delta_{n+m}}
  {\delta_{n}\delta_{m}}\hat{f}(n)\hat{g}(m)z^{n+m}.
  \end{eqnarray}

Note that $\ell^{p}(\beta)$ is not unitarily equivalent to
$\ell^{p}(\widetilde{\beta})$ with $\widetilde{\beta}(n):=
\delta_{n+1}\beta(n)/(\delta_{1}\delta_{n})$, because in general
the sequence $\{\widetilde{\beta}(n)/\beta(n)\}$ is not
necessarily constant (see \cite{sh}). Let $u, \varphi\in l^0(\beta)$, the
set of all formal power series. The operator $u\diamond
C_\varphi:\ell^{p}(\beta)\rightarrow l^0(\beta)$ defined by
$$(u\diamond C_\varphi)(f)=u\diamond (f\circ\varphi), \ \   \ \ \ \ \ f\in \ell^{p}(\beta)$$
is called $\diamond$-substitution operator on $\ell^{p}(\beta)$.
Note that if
$\lim_{n}{\beta(n+1)}/{\beta(n)}=1$ or
$\liminf_{n}\beta(n)^{{1}/{n}}=1$, we shall assume that
$\varphi(\mathbb{D})\subseteq\mathbb{D}$, where $\mathbb{D}$ is
the open unit disc.

\vspace*{0.3cm} Substitution (weighted composition) operators are
well-studied in many classical spaces. In this note we initiate
the study of $\diamond$-substitution operators on
$\ell^{p}(\beta)$ with $\diamond$-multiplicaton. There is a few
results about substitution operator $uC_\varphi$  on
$\ell^{p}(\beta)$ (see \cite{g,z}). In the next section, by an
inductive argument on the pair $(u, \varphi)$, we give some
sufficient and necessary conditions for boundedness of
$\diamond$-substitution operator $u\diamond C_\varphi$ on
$\ell^{p}(\beta)$.

\section{\sc\bf Boundedness of $u\diamond
C_\varphi$ on $\ell^{p}(\beta)$}

\vspace*{0.3cm} {\bf Theorem 2.1.} Let $m, n\in\mathbb{N}_0$.
If $\varphi(z)=z^{m}$, then $C_{\varphi}\in B(l^{p}(\beta))$ if and
only if $\sup_{n\in\mathbb{N}_0}{\beta(nm)}/{\beta(n)}<\infty$. In this
case
$$\|C_{\varphi}\|=\sup_{n\in\mathbb{N}_0}\frac{\beta(nm)}{\beta(n)} .$$

\vspace*{0.3cm} {\bf Proof.} Suppose that
$M:=\sup_{n\in\mathbb{N}_0}{\beta(nm)}/{\beta(n)}$ is finite. Then for every
$f\in l^{p}(\beta)$ we have
$$f\circ\varphi(z)=f(\varphi(z))=f(z^m)=\sum^{\infty}_{n=0}\hat{f}(n)z^{nm}$$
and
$$\|C_{\varphi}f\|^{p}_{\beta}=\|\sum^{\infty}_{n=0}\hat{f}(n)z^{nm}\|^{p}_{\beta}=
\sum^{\infty}_{n=0}|\hat{f}(n)|^{p}\beta(nm)^{p}$$
$$
\leq M^{p}\sum^{\infty}_{n=0}|\hat{f}(n)|^{p}\beta(n)^{p}=
M^{p}\|f\|^{p}_{\beta}.$$
Thus $\|C_{\varphi}\|\leq M$,
and so $C_{\varphi}$ is a bounded operator on $l^{p}(\beta)$.
On the other hand, if we put
$f_{n}(z)=z^{n}$, then we get that
$$\beta(nm)^{p}=\|z^{nm}\|^{p}_{\beta}=
 \|C_{\varphi}z^{n}\|^{p}_{\beta}\leq\|C_{\varphi}\|^{p}\|z^{n}\|^{p}_{\beta}=\|C_{\varphi}\|^{p}\beta(n)^{p},$$
and so $M\leq\|C_{\varphi}\| $. \ $\Box$

 \vspace*{0.3cm} {\bf Theorem 2.2.} Let $\alpha_{i}\in\mathbb{R}$ and
$\varphi(z)=\sum^{d}_{m=0}\alpha_{m}z^{m}$. Put
$$B_{d}=\sum^{\infty}_{n=0}\left[\sum_{
L\in\mathbb{N}_0}\left(\frac{|\theta_{n,L}|\beta(n))}{\beta(L)}\right)^{q}\right]^{\frac{p}{q}},$$
where $q$ is the conjugate exponent to $p$,
$$\theta_{n,L}:=\sum_{(l_{0},l_{1}, ...,l_{d})\in
N(n,L)}\frac{L!\alpha^{l_{0}}_{0}\alpha^{l_{1}}_{1}...\alpha^{l_{d}}_{d}}{(l_{0})_{!}...(l_{d})_{!}} ,$$
with $N(n,L):=\{(l_{0},l_{1}, ...,l_{d})| \ \
l_{i}\in\mathbb{N}_0, \ \ \sum^{d}_{i=0}l_{i}=L,\ \
\sum^{d}_{i=1}il_{i}=n\}$.
Then the following hold.\\

(a) If $B_{d}<\infty$, then $C_{\varphi}\in B(l^{p}(\beta))$ and $\|C_{\varphi}\|\leq \sqrt[p]{B_{d}}$.

(b) If $C_{\varphi}\in B(l^{p}(\beta))$, then
$$M:=\sup_{n\geq0}\frac{\left(\sum_{L\in\mathbb{N}_0}|\theta_{n,L}|^{p}\beta(L)^p\right)^{\frac{1}{p}}}{\beta(n)}<\infty ,$$
and $\|C_{\varphi}\|\leq M$.

\vspace*{0.3cm} {\bf Proof.} (a) Take
$N(n)=\{(l_{0},l_{1},...,l_{d})| \ \ l_{i}\in\mathbb{N}_0, \ \
\Sigma^{d}_{i=0}l_{i}=n\}$.
 Then for each $f\in l^{p}(\beta)$, we have
$$f\circ\varphi(z)=f(\varphi(z))=\sum^{\infty}_{n=0}\hat{f}(n)(\varphi(z))^{n}=\sum^{\infty}_{n=0}\hat{f}(n)(\sum^{d}_{m=0}\alpha_{m}z^{m})^{n}$$
$$=
\sum^{\infty}_{n=0}\hat{f}(n)\left(\sum_{(l_{0}, l_{1},...,l_{d})\in
N(n)}\frac{n!\alpha^{l_{0}}_{0}\alpha^{l_{1}}_{1}...\alpha^{l_{d}}_{d}}{(l_{0})_{!}...(l_{d})_{!}}
z^{\sum^{d}_{i=1}il_{i}}\right).$$
This implies that
$$\widehat{f\circ\varphi}(n)=\sum_{ L\in\mathbb{N}_0}\hat{f}(L)\left(\sum_{(l_{0},l_{1}, ...,l_{d})\in
N(n,L)}\frac{L!\alpha^{l_{0}}_{0}\alpha^{l_{1}}_{1}...\alpha^{l_{d}}_{d}}{(l_{0})_{!}...(l_{d})_{!}}
\right)=\sum_{ L\in\mathbb{N}_0}\hat{f}(L)\theta_{n,L},$$
and
$$\|C_{\varphi}f\|^{p}_{\beta}=\sum^{\infty}_{n=0}|\widehat{f\circ\varphi}(n)|^p|\beta(n)^{p}=
\sum^{\infty}_{n=0}|\sum_{
L\in\mathbb{N}_0}\hat{f}(L)\theta_{n,L}|^{p}\beta(n)^{p}.$$
Now, by using H\"{o}lder inequality we get that
$$\|C_{\varphi}f\|^{p}_{\beta}\leq\sum^{\infty}_{n=0}\left(\sum_{ L\in\mathbb{N}_0}
|\hat{f}(L)|\beta(L)\frac{|\theta_{n,L}|\beta(n)}{\beta(L)}\right)^{p}
$$
$$\leq\sum^{\infty}_{n=0}\left(\sum_{ L\in\mathbb{N}_0}|\hat{f}(L)|^{p}\beta(L)^{p}\right)\left(\sum_{
L\in\mathbb{N}_0}(\frac{|\theta_{n,L}|\beta(n))}{\beta(L)})^{q}\right)^{\frac{p}{q}}.$$
It follows that $\|C_{\varphi}f\|_{\beta}\leq \sqrt[p]{B_{d}}\ \|f\|_{\beta}$,
and so $\|C_{\varphi}\|\leq\sqrt[p]{B_{d}}$.

\vspace*{0.3cm}
(b) Let $C_{\varphi}\in B(l^{p}(\beta))$. For $f_{n}(z)=z^n$ we have
$$f_{n}(\varphi(z))=(\sum^{d}_{m=0}\alpha_{m}z^{m})^n=\sum_{(l_{0}, l_{1},...,l_{d})\in
N(n)}\frac{n!\alpha^{l_{0}}_{0}\alpha^{l_{1}}_{1}...\alpha^{l_{d}}_{d}}{(l_{0})_{!}...(l_{d})_{!}}
z^{\sum^{d}_{i=1}il_{i}}.$$
 So, for each $L\in \mathbb{N}_0$
$$\widehat{f_{n}\circ\varphi}(L)=\sum_{(l_{0}, l_{1},...,l_{d})\in
N(L,n)}\frac{n!\alpha^{l_{0}}_{0}\alpha^{l_{1}}_{1}...\alpha^{l_{d}}_{d}}{(l_{0})_{!}...(l_{d})_{!}}=\theta_{n,L}.$$
Thus we have
$$\|C_{\varphi}f_{n}\|^p_{\beta}=\sum_{L\in\mathbb{N}_0}|\theta_{n,L}|^p\beta(L)^p\leq\|C_{\varphi}\|^p\|f_{n}\|^p_{\beta}=\|C_{\varphi}\|^p\beta(n)^p,$$
and so
$M\leq\|C_{\varphi}\|<\infty$.

\vspace*{0.3cm} {\bf Theorem 2.3.} Let $m\in\mathbb{N}_0$,
$M(m)=\{km: \ k\in\mathbb{N}_0\}$ and let $\varphi(z)=z^{m}$.

\vspace*{0.3cm} (a) Put $$\alpha=\sup_{n\in\mathbb{N}_0}\ \sum_{k\in A_{m,
n}}\left(\frac{\delta_{n}\beta(n)}{\delta_{k}\delta_{n-k}\beta(k)\beta(\frac{n-k}{m})}\right)^{q},$$
where $A_{n,m}:=\{k: 0\leq k\leq n \ \mbox {and} \  n-k \in
M(m)\}$ and $1/p+1/q=1$. If $u\in \ell^p(\beta)$ and
$\alpha<\infty$, then $u\diamond C_{\varphi}$ is a bounded
operator on $\l^{p}(\beta)$. In this case $\|u\diamond
C_{\varphi}\|\leq\alpha^{{1}/{q}}\|u\|_{\beta}$.

\vspace*{0.3cm} (b) If $u\diamond C_{\varphi}$ is a bounded operator
on $\l^{p}(\beta)$, then
$$M:=\sup_{l\in\mathbb{N}_0 }\frac{\left(\sum^{\infty}_{n=ml}\left(\frac{\delta_{n}\beta(n)|\hat{u}(n-ml)|}
{\delta_{ml}\delta_{n-ml}}\right)^p\right)^{\frac{1}{p}}}{\beta(l)}<\infty,$$
and $M\leq\|u\diamond C_{\varphi}\|$.

\vspace*{0.3cm} {\bf Proof.} (a) Let $u, f\in l^{p}(\beta)$ and
$\alpha<\infty$. Then
$f\circ\varphi(z)=\sum^{\infty}_{n=0}\hat{f}(n)z^{nm}$,
$\widehat{f\circ\varphi}(n-k)=\hat{f}(\frac{n-k}{m})$ and
$$\widehat{u\diamond (f\circ\varphi)}(n)=\sum^{n}_{k=0}\frac{\delta_{n}\hat{u}(k)\widehat{f\circ\varphi}(n-k)}
{\delta_{k}\delta_{n-k}}=\sum_{k\in
A_{n,m}}\frac{\delta_{n}\hat{u}(k)\hat{f}(\frac{n-k}{m})}{\delta_{k}\delta_{n-k}}.$$
It follows that
$$\|u\diamond
 (f\circ\varphi)\|^{p}_{\beta}=\sum^{\infty}_{n=0}\mid\widehat{u\diamond
(f\circ\varphi)}(n)\beta(n)\mid^{p}=\sum^{\infty}_{n=0}|\sum_{k\in
A_{n,m}}\frac{\delta_{n}\beta(n)}{\delta_{k}\delta_{n-k}}\hat{u}(k)\hat{f}(\frac
{n-k}{m})|^{p}$$
$$\leq\sum^{\infty}_{n=0}\left(\sum_{k\in A_{n, m}}\frac{\delta_{n}\beta(n)}{\delta_{k}\delta_{n-k}\beta(k)\beta(\frac{n-k}{m})}|\hat{u}(k)|\beta(k)|\hat{f}(\frac
{n-k}{m})|\beta(\frac{n-k}{m})\right)^{p}$$
$$\overset{{\rm H\ddot{o}lder}}{\leq}\sum^{\infty}_{n=0}\left(\sum_{k\in A_{n, m}}\left(\frac{\delta_{n}\beta(n)}{\delta_{k}\delta_{n-k}\beta(k)\beta(\frac{n-k}{m})}\right)^{q}
\right)^{\frac{p}{q}}$$$$\times\left(\sum_{k\in A_{n, m}}|\hat{u}(k)|^{p}\beta(k)^p|\hat{f}(\frac
{n-k}{m})|^{p}\beta(\frac{n-k}{m})^{p}\right)$$
$$\leq\alpha^{\frac{p}{q}}\sum^{\infty}_{n=0}\sum_{k\in A_{n, m}}|\hat{u}(k)|^{p}\beta(k)^p|\hat{f}(\frac
{n-k}{m})|^{p}\beta(\frac{n-k}{m})^{p}$$
$$\leq\alpha^{\frac{p}{q}}\|f\|^{p}_{\beta}\|u\|^p_{\beta}.$$
Thus $\|u\diamond C_{\varphi}\|\leq\alpha^{{1}/{q}}\|u\|_{\beta}$,
and so $u\diamond C_{\varphi}$ is bounded.

\vspace*{0.3cm} (b) Let $u\diamond C_{\varphi}\in B(l^{p}(\beta))$. For $f_{l}(z)=z^l$ we have
$f_{l}\circ\varphi(z)=z^{ml}$ and
$$\widehat{u\diamond f_{l}\circ\varphi}(n)=\sum^n_{k=0}\frac{\delta_{n}\hat{u}(k)\hat{f_{l}}(\frac{n-k}{m})}{\delta_{k}\delta_{n-k}}.$$
It follows that
\[\widehat{u\diamond f_{l}\circ\varphi}(n) = \left \{ \begin{array}{ll}
\frac{\delta_{n}\hat{u}(n-ml)}{\delta_{ml}\delta_{n-ml}} &      \ n\geq ml\\
 0      & \ n<ml.
\end{array}
\right. \]
Thus
$$\|\widehat{u\diamond f_{l}\circ\varphi}\|^p_{\beta}=\sum^{\infty}_{n=ml}\left(\frac{\delta_{n}\beta(n)|\hat{u}(n-ml)|}{\delta_{ml}\delta_{n-ml}}\right)^p
$$$$\leq\|u\diamond C_{\varphi}\|^p\|f_{l}\|^p_{\beta}=\|u\diamond C_{\varphi}\|^p\beta(l)^p.$$
This implies that
$$\frac{\left(\sum^{\infty}_{n=ml}\left(\frac{\delta_{n}\beta(n)|\hat{u}(n-ml)|}
{\delta_{ml}\delta_{n-ml}}\right)^p\right)^{\frac{1}{p}}}{\beta(l)}\leq\|u\diamond C_{\varphi}\|<\infty,$$
and so $M<\|u\diamond C_{\varphi}\|<\infty$.\ $\Box$

\vspace*{0.3cm}
For each $u\in\ell^{p}(\beta)$, let $M_{\diamond,
u}: \ell^{p}(\beta)\rightarrow \ell^{0}(\beta)$ defined by
$M_{\diamond, u}(f)=u\diamond f$ be its corresponding
$\diamond$-multiplication linear operator. Put
$$\alpha_o:=\sup_{n\in\mathbb{N}_0}\sum^{n}_{k=0}\left(\frac{\delta_{n}\beta(n)}
{\delta_{k}\delta_{n-k}\beta(k)\beta(n-k)}\right)^{q}.$$
Then by Theorem 2.3, if $\alpha_0<\infty$ then
$\|M_{\diamond, u}\|\leq{\alpha_o}^{{1}/{q}}$.
Also, it is easy to see that the constant
function $f=1$ is a unity for $(\ell^{p}(\beta),\diamond)$.
These observations establish the following corollary.

\vspace*{0.3cm} {\bf Corollary 2.4.}
$(\ell^{p}(\beta),\diamond)$
is a unital commutative Banach algebra.

\vspace*{0.3cm} {\bf Theorem 2.5.} Let $d\in\mathbb{N}_0$,
$\alpha_{i}\in\mathbb{R}$, $u(z)=z^{m_{0}}$ and
$\varphi(z)=\sum^{d}_{m=0}\alpha_{m}z^{m}$.

\vspace*{0.3cm} (a) If
$$\alpha:=\sum^{\infty}_{n=m_{0}}\left(\frac{\delta_{n}\beta(n)}{\delta_{m_{0}}\delta_{n-m_{0}}}\right)^p\left(\sum_{L\in\mathbb{N}_0}
(\frac{|\theta_{n-m_{0},L}|}{\beta(L)})^q\right)^{\frac{p}{q}}<\infty,$$
then $u\diamond
C_{\varphi}$ is a bounded operator on $\l^{p}(\beta)$.

\vspace*{0.3cm} (b) If $u\diamond C_{\varphi}$ is a bounded operator
on $\l^{p}(\beta)$, then
$$M:=\sup_{l\geq0}\frac{\left[\sum^{\infty}_{n=m_{0}}\left(\frac{\delta_{n}|\theta_{n-m_{0},l}|\beta(n)}
{\delta_{m_{0}}\delta_{n-m_{0}}}\right)^p\right]^{\frac{1}{p}}}{\beta(l)}<\infty$$
with $M\leq\|u\diamond C_{\varphi}\|$.

\vspace*{0.3cm} {\bf Proof.} (a) Let $f\in l^{p}(\beta)$. Then we have
$$f\circ\varphi(z)=f(\varphi(z))=\sum^{\infty}_{n=0}\hat{f}(n)(\varphi(z))^{n}=\sum^{\infty}_{n=0}\hat{f}(n)(\sum^{d}_{m=0}\alpha_{m}z^{m})^{n}$$
$$=
\sum^{\infty}_{n=0}\hat{f}(n)\left(\sum_{(l_{0}, l_{1},...,l_{d})\in
N(n)}\frac{n!\alpha^{l_{0}}_{0}\alpha^{l_{1}}_{1}...\alpha^{l_{d}}_{d}}{(l_{0})_{!}...(l_{d})_{!}}
z^{\sum^{d}_{i=1}il_{i}}\right).$$
This implies that
$$\widehat{f\circ\varphi}(n)=\sum_{ L\in\mathbb{N}_0}\hat{f}(L)\left(\sum_{(l_{0},l_{1}, ...,l_{d})\in
N(n,L)}\frac{L!\alpha^{l_{0}}_{0}\alpha^{l_{1}}_{1}...\alpha^{l_{d}}_{d}}{(l_{0})_{!}...(l_{d})_{!}}\right)=\sum_{
L\in\mathbb{N}_0}\hat{f}(L)\theta_{n,L},$$
and
 $$\widehat{u\diamond
 (f\circ\varphi)}(n)=\sum^{n}_{k=0}
\frac{\delta_{n}}{\delta_{k}\delta_{n-k}}\hat{u}(k)\widehat{f\circ\varphi}(n-k).$$
Since $u(z)=z^{m_{0}}$, we get that
\[\widehat{u\diamond (f\circ\varphi)}(n) = \left \{ \begin{array}{ll}
\frac{\delta_{n}}{\delta_{m_{0}}\delta_{n-m_{0}}}\sum_{
L\in\mathbb{N}_0}\hat{f}(L)\theta_{n-m_{0},L} &      \ n\geq m_{0}\\
 0      & \ n<m_{0}.
\end{array}
\right. \]
Thus
$$\|u\diamond
 (f\circ\varphi)\|^{p}_{\beta}=\sum^{\infty}_{n-m_{0}}\frac{\delta^p_{n}}{\delta^p_{m_{0}}\delta^p_{n-m_{0}}}
\mid\sum_{
L\in\mathbb{N}_0}\hat{f}(L)\theta_{n-m_{0},L}\mid^p\beta(n)^p$$

$$=\sum^{\infty}_{n=m_{0}}\frac{\delta^p_{n}\beta(n)^p}{\delta^p_{m_{0}}\delta^p_{n-m_{0}}}
|\sum_{
L\in\mathbb{N}_0}\hat{f}(L)\beta(L)\frac{\theta_{n-m_{0},L}}{\beta(L)}|^p.$$
Now, by H\"{o}lder inequality we have
$$\|u\diamond
(f\circ\varphi)\|^p_{\beta}\leq\|f\|^p_{\beta}\sum^{\infty}_{n=m_{0}}\frac{\delta^p_{n}
\beta(n)^p}{\delta^p_{m_{0}}\delta^p_{n-m_{0}}}\left(\sum_{L\in\mathbb{N}}
(\frac{|\theta_{n-m_{0},L}|}{\beta(L)})^q\right)^{\frac{p}{q}}.$$
Thus $\|u\diamond
C_{\varphi}\|^p\leq\alpha$, and so $u\diamond C_{\varphi}$ is bounded.

\vspace*{0.3cm}
(b) Let $u\diamond C_{\varphi}\in B(l^{p}(\beta))$. Let $f_{l}(z)=z^l$,
then by part (a) we have
$$\widehat{u\diamond
(f_{l}\circ\varphi)}(n)=
\frac{\delta_{n}}{\delta_{m_{0}}\delta_{n-m_{0}}}\sum_{
L\in\mathbb{N}_0}\hat{f_{l}}(L)\left(\sum_{(l_{0},l_{1},
...,l_{d})\in
N(n-m_{0},L)}\frac{L!\alpha^{l_{0}}_{0}\alpha^{l_{1}}_{1}...\alpha^{l_{d}}_{d}}{(l_{0})_{!}...(l_{d})_{!}}\right)$$$$=
\frac{\delta_{n}}{\delta_{m_{0}}\delta_{n-m_{0}}}\left(\sum_{(l_{0},l_{1},
...,l_{d})\in
N(n-m_{0},L)}\frac{L!\alpha^{l_{0}}_{0}\alpha^{l_{1}}_{1}...\alpha^{l_{d}}_{d}}{(l_{0})_{!}...(l_{d})_{!}}\right).$$
Hence
$$\|u\diamond
(f_{l}\circ\varphi)\|^{p}_{\beta}=\sum^{\infty}_{n=m_{0}}\left[\frac{\delta_{n}|\theta_{n-m_{0},L}|\beta(n)}
{\delta_{m_{0}}\delta_{n-m_{0}}}\right]^p\leq\|u\diamond
C_{\varphi}\|^p\beta(l)^p,$$
so
$M\leq\|u\diamond
C_{\varphi}\|$.
This completes proof.

\vspace*{0.3cm} {\bf Corollary 2.6.} Let $m_{1},
m_{2}\in\mathbb{N}_0$, $u(z)=z^{m_{1}}$, and let
$\varphi(z)=z^{m_{2}}$.

\vspace*{0.3cm} (a) Put
$$\gamma=\sup_{n\in B_{m_1, m_2}}\frac{\delta_{n}\beta(n)}{\delta_{m_{1}}\delta_{n-m_{1}}\beta(\frac{n-m_{1}}{m_{2}})} \ ,$$
where $B_{m_1, m_2}:=\{n: n\geq m_{1}\ {\rm and}\  n-m_{1} \in
M(m_{2})\}$. If $\gamma<\infty$, then $u\diamond C_{\varphi}$ is a
bounded operator on $l^{p}(\beta)$ and $\|u\diamond
C_{\varphi}\|\leq\gamma$.

\vspace*{0.3cm} (b) If $u\diamond C_{\varphi}$ is a bounded operator
on $\l^{p}(\beta)$, then
$$K:=\sup_{m\in\mathbb{N}_0}\frac{\delta_{m_{1}+mm_{2}}\beta(m_{1}+mm_{2})}{\delta_{mm_{2}}\delta_{m_{1}}\beta(m)}<\infty,$$
with $K\leq\|u\diamond C_{\varphi}\|$.

\vspace*{0.3cm} {\bf Proof.} (a) Let $f\in l^{p}(\beta)$. Then precisely the same calculation as
that shown in the proof of Theorem 2.5 yields that
\[\widehat{u\diamond (f\circ\varphi)}(n) = \left \{ \begin{array}{ll}
\frac{\delta_{n}\hat{f}(\frac{n-m_{1}}{m_{2}})}{\delta_{m_{1}}\delta_{n-m_{1}}} &      \ n\in B_{m_{1},m_{2}}\\
 0      & \ n\notin B_{m_{1},m_{2}}.
\end{array}
\right. \]
Then
$$\|u\diamond
(f\circ\varphi)\|^{p}_{\beta}=\sum_{n\in B_{m_1,
m_2}}(\frac{\delta_{n}\beta(n)}{\delta_{m_{1}}\delta_{n-m_{1}}})^{p}|\hat{f}(\frac{n-m_{1}}{m_{2}})|^{p}$$
$$=\sum_{n\in B_{m_1, m_2}}\left(\frac{\delta_{n}\beta(n)}{\delta_{m_{1}}\delta_{n-m_{1}}\beta(\frac{n-m_{1}}{m_{2}})}\right)^{p}|\hat{f}
(\frac{n-m_{1}}{m_{2}})|^{p}\beta(\frac{n-m_{1}}{m_{2}})^{p}$$
$$\leq\gamma^p\sum_{n\in B_{m_1, m_2}}|\hat{f}
(\frac{n-m_{1}}{m_{2}})|^{p}\beta(\frac{n-m_{1}}{m_{2}})^{p}\leq\gamma^p\|f\|^{p}_{\beta}.$$
This implies that $\|u\diamond C_{\varphi}\|\leq\gamma$.

\vspace*{0.3cm} (b) Let $u\diamond C_{\varphi}\in
B(l^{p}(\beta))$. Let $f_{m}(z)=z^m$, then we get that
\[\widehat{u\diamond (f_{m}\circ\varphi)}(n) = \left \{ \begin{array}{ll}
\frac{\delta_{mm_{2}}+m_{1}}{\delta_{m_{1}}\delta_{mm_{2}}} &      \ n=mm_{2}+m_{1}\\
 0      & \ \mbox {otherwise}.
\end{array}
\right. \]
 Thus
 $$\|u\diamond
(f_{m}\circ\varphi)\|^{p}_{\beta}=\left(\frac{\delta_{mm_{2}+m_{1}}\beta(mm_{2}+m_{1})}{\delta_{m_{1}}\delta_{mm_{2}}}\right)^p\leq\|u\diamond
C_{\varphi}\|^p \beta(m)^p.$$
This implies that
$K\leq\|u\diamond C_{\varphi}\|<\infty$.\ $\Box$



\vspace*{0.3cm}


\begin{thebibliography}{99}

\bibitem{c}
C. Cowen and Barbara MacCluer, Composition operators on spaces of
analytic functions, CRC Press, Boca Raton, FL, 1995.

\bibitem{g}
Gajath K. Gunatillake,  Weighted composition operators, Ph.D.
Thesis, 2005.


\bibitem{sh}
A. L. Shields, Weighted shift operators and analytic function
theory, Math. Surveys, A. M. S. Providence, {\bf 13} (1974),
49-128.

\bibitem{z}
N. Zorboska, Composition operators on weighted Hardy spaces,
Ph.D. Thesis, 1988.

\end{thebibliography}
\end{document}